\documentclass{gtart}

\def\ifplaintex{\expandafter\ifx\csname documentclass\endcsname\relax}

\def\gtp{{\mathsurround=0pt\it $\cal G\mskip-2mu$eometry \&\ 
$\cal T\!\!$opology $\cal P\!$ublications}}  % GT publications

\def\recd{{\small Received:\qua\receiveddate\ifx\reviseddate\relax
\else\qquad Revised:\qua\reviseddate\fi\par}} 

\def\lognumber#1{\def\thelognumber{#1}}
\def\volumenumber#1{\def\thevolumenumber{#1}}
\def\volumeyear#1{\def\thevolumeyear{#1}}
\def\papernumber#1{\def\thepapernumber{#1}}
\def\pagenumbers#1#2{\def\startpage{#1}\def\finishpage{#2}}
\def\published#1{\def\publishdate{#1}}

\def\received#1{\def\receiveddate{#1}}

\def\accepted#1{\def\accepteddate{#1}}

\let\\\par\let\thelognumber\relax\let\thevolumenumber\relax
\let\thepapernumber\relax\let\thevolumeyear\relax\let\startpage\relax
\let\finishpage\relax\let\publishdate\relax\let\receiveddate\relax
\let\reviseddate\relax\let\accepteddate\relax\let\theasciititle\relax
\let\theasciiauthors\relax
\let\theasciiabstract\relax

\let\theasciiemail\relax

\ifplaintex
\font\logobig=cmssbx10 scaled 3836
\font\logomed=cmssbx10 scaled 2557
\else
\font\logobig=cmssbx10 scaled 4200
\font\logomed=cmssbx10 scaled 2800
\fi

\long\def\makeagttitle{   %%% start of definition of \makeagttitle
\count0=\startpage
\agt\hfill      %   Journal title (top left) 
\hbox to 45truept{\vbox to 0pt{\vglue -13truept{\logomed A\kern -.37em{\logobig 
T}\kern -.38em G}\vss}\hss}
\break
{\small Volume \thevolumenumber\ (\thevolumeyear)
\startpage--\finishpage\nl
Published: \publishdate}

\vglue .25truein

{\parskip=0pt\leftskip 0pt plus
1fil\def\\{\par\smallskip}{\Large\bf\thetitle}\par\medskip} \vglue
0.05truein

{\parskip=0pt\leftskip 0pt plus 1fil\def\\{\par}{\sc\theauthors}
\par\medskip}%
 
\vglue 0.03truein 

{\small\leftskip 25truept\rightskip 25truept{\bf Abstract}\stdspace\theabstract

{\bf AMS Classification}\stdspace\theprimaryclass
\ifx\thesecondaryclass\relax\else; \thesecondaryclass\fi\par
{\bf Keywords}\stdspace \thekeywords\par}\vglue 7truept

}   %%%% end of definition of \makeagttitle

\ifplaintex
\font\phead=cmsl9 scaled 950
\font\pnum=cmbx10 scaled 913
\font\pfoot=cmsl9 scaled 950
\headline{\vbox to 0pt{\vskip -4.5mm\line{\small\phead\ifnum
\count0=\startpage ISSN 1472-2739 (on-line) 1472-2747 (printed)
\hfill {\pnum\folio}\else\ifodd\count0\def\\{ }% 
\ifx\theshorttitle\relax\thetitle\else\theshorttitle\fi\hfill{\pnum\folio}
\else\def\\{ and }{\pnum\folio}\hfill\ifx\theshortauthors\relax\theauthors
\else\theshortauthors\fi\fi\fi}\vss}}
\footline{\vbox to 0pt{\vglue 0mm\line{\small\pfoot\ifnum\count0=\startpage
\copyright\ \gtp\hfill\else
\agt, Volume \thevolumenumber\ (\thevolumeyear)\hfill\fi}\vss}}
\else
\font=cmsl9 scaled 1050
\font\lnum=cmbx10 
\font=cmsl9 scaled 1050
\makeatletter
\def\@oddhead{{\small\lhead\ifnum\count0=\startpage ISSN 1472-2739 
(on-line) 1472-2747 (printed)\hfill {\lnum\number\count0}\else\ifodd\count0
\def\\{ }\ifx\theshorttitle\relax \thetitle \else\theshorttitle\fi\hfill
{\lnum\number\count0}\else\def\\{ and }{\lnum\number\count0}
\hfill\ifx\theshortauthors\relax 
\theauthors\else\theshortauthors\fi\fi\fi}}\def\@evenhead{\@oddhead}
\def\@oddfoot{\small\lfoot\ifnum\count0=\startpage\copyright\ \gtp\hfill\else
\agt, Volume \thevolumenumber\ (\thevolumeyear)\hfill\fi}
\def\@evenfoot{\@oddfoot}
\makeatother
\fi
\let\maketitlepage\makeagttitle
\let\makeshorttitle\maketitlepage
\let\maketitle\maketitlepage

\newwrite\gtoutfile
\long\gdef\makeheadfile{  %%% start of definition of \makeheadfile
{\def\\{, }\def\s{ }
\immediate\openout\gtoutfile head.xxx
\immediate\write\gtoutfile{To: math@arxiv.org}
\immediate\write\gtoutfile{Subject: put OR rep NNNNN:ppppp}
\immediate\write\gtoutfile{--text follows this line--}
\immediate\write\gtoutfile{Proxy-for: \ifx\theasciiauthors\relax
\theauthors\else\theasciiauthors\fi\s<\ifx\theasciiemail\relax\theemail\else\theasciiemail\fi>}
\immediate\write\gtoutfile{\noexpand\\}
\immediate\write\gtoutfile{Authors: \ifx\theasciiauthors\relax
\theauthors\else\theasciiauthors\fi}
{\def\\{ }\immediate\write\gtoutfile{Title: \ifx\theasciititle\relax
\thetitle\else\theasciititle\fi}}
\immediate\write\gtoutfile{Subj-class: GT or SG, GR etc}
\immediate\write\gtoutfile{MSC-class: \theprimaryclass\ifx\thesecondaryclass\relax\else, \thesecondaryclass\fi}
\immediate\write\gtoutfile{Journal-ref: Algebr. Geom. Topol. \thevolumenumber\s
(\thevolumeyear) \startpage-\finishpage}
\immediate\write\gtoutfile{Comments: Published by Algebraic and
Geometric Topology at}
\immediate\write\gtoutfile{\s\s\s  http://www.maths.warwick.ac.uk/agt/AGTVol\thevolumenumber/agt-\thevolumenumber-\thepapernumber.abs.html}
\immediate\write\gtoutfile{\noexpand\\}
\immediate\write\gtoutfile{}
\ifx\theasciiabstract\relax
\immediate\write\gtoutfile{\theabstract}\else
\immediate\write\gtoutfile{\theasciiabstract}\fi
\immediate\write\gtoutfile{}
\immediate\write\gtoutfile{\noexpand\\}
\immediate\write\gtoutfile{}
\immediate\closeout\gtoutfile}}  %%% end of definition of \makeheadfile

\def\maketitlepage{\makeagttitle\makeheadfile}
\let\makeshorttitle\maketitlepage
\let\maketitle\maketitlepage

\lognumber{8}
\volumenumber{3}
\volumeyear{2003}
\papernumber{8}
\published{22 February 2003}
\pagenumbers{187}{205}
\received{9 January 2003}
\accepted{13 February 2003}

\usepackage{amsmath,amssymb}

\input xypic

\newcommand{\Alg}       {\operatorname{Alg}}
\newcommand{\EPR}       {\operatorname{EPR}}

\newcommand{\FG}        {\operatorname{FG}}
\newcommand{\GFG}       {\operatorname{GFG}}
\newcommand{\LEPR}      {\operatorname{LEPR}}
\newcommand{\LFG}       {\operatorname{LFG}}
\newcommand{\DMod}      {\operatorname{DMod}}
\newcommand{\EPA}       {\operatorname{EPA}}
\newcommand{\Mod}       {\operatorname{Mod}}
\newcommand{\Grass}     {\operatorname{Grass}}
\newcommand{\Thom}      {\operatorname{Thom}}

\newcommand{\ext}       {{\operatorname{ext}}}
\newcommand{\op}        {{\operatorname{op}}}
\newcommand{\spec}      {\operatorname{spec}}
\newcommand{\spf}       {\operatorname{spf}}

\newcommand{\Gm}        {\Gamma}
\newcommand{\al}        {\alpha}
\newcommand{\bt}        {\beta}
\newcommand{\gm}        {\gamma}
\newcommand{\lm}        {\lambda}
\newcommand{\sg}        {\sigma}
\newcommand{\om}        {\omega}
\newcommand{\tht}       {\theta}

\newcommand{\CB}        {{\mathcal{B}}}
\newcommand{\CE}        {{\mathcal{E}}}
\newcommand{\CLC}       {{\mathcal{L}_\mathbb{C}}}
\newcommand{\CU}        {{\mathcal{U}}}
\newcommand{\CV}        {{\mathcal{V}}}
\newcommand{\CR}        {{\mathcal{R}}}

\newcommand{\C}         {{\mathbb{C}}}
\newcommand{\CP}        {\mathbb{C}P}
\newcommand{\CPi}       {\mathbb{C}P^\infty}

\newcommand{\OG}        {{\mathcal{O}}_G}
\newcommand{\OS}        {{\mathcal{O}}_S}
\newcommand{\Sg}        {\Sigma}
\newcommand{\Fp}        {{\mathbb{F}_p}}          
\newcommand{\N}         {{\mathbb{N}}}
\newcommand{\Z}         {{\mathbb{Z}}}
\newcommand{\psb}[1]    {[\![#1]\!]}
\newcommand{\sixth}     {{\textstyle\frac{1}{6}}}
\newcommand{\sm}        {\setminus}
\newcommand{\st}        {\;|\;}
\newcommand{\xra}       {\xrightarrow}

\newcommand{\ot}        {\otimes}
\newcommand{\tE}        {\widetilde{E}}

\newcommand{\tx}        {\widetilde{x}}
\newcommand{\tm}        {\times}
\newcommand{\xla}       {\xleftarrow}
\newcommand{\Smash}     {\wedge}
\newcommand{\bigWedge}  {\bigvee}
\newcommand{\sse}       {\subseteq}
\newcommand{\ths}       {\ltimes}
\newcommand{\tp}        {\tilde{p}}
\newcommand{\tq}        {\tilde{q}}
\newcommand{\tr}        {\tilde{r}}
\newcommand{\ts}        {\tilde{s}}
\newcommand{\ux}        {\underline{x}}
\newcommand{\half}      {{\textstyle\frac{1}{2}}}

\newcommand{\zin}       {1\Smash\eta}
\newcommand{\zni}       {\eta\Smash 1}
\newcommand{\znn}       {\eta\Smash\eta}
\newcommand{\ziin}      {1\Smash 1\Smash\eta}
\newcommand{\zini}      {1\Smash\eta\Smash 1}
\newcommand{\znii}      {\eta\Smash 1\Smash 1}
\newcommand{\znni}      {\eta\Smash\eta\Smash 1}
\newcommand{\znin}      {\eta\Smash 1\Smash\eta}
\newcommand{\zinn}      {1\Smash\eta\Smash\eta}

\newcommand{\hocolim}{\operatornamewithlimits{\underset{\longrightarrow}{holim}}}

\renewcommand{\:}{\colon}

\newtheorem{theorem}{Theorem}[section]

\newtheorem{lemma}[theorem]{Lemma}
\newtheorem{proposition}[theorem]{Proposition}
\newtheorem{corollary}[theorem]{Corollary}
\theoremstyle{definition}
\newtheorem{assumption}[theorem]{Assumption}
\newtheorem{remark}[theorem]{Remark}
\newtheorem{definition}[theorem]{Definition}

\begin{document}
\title{Realising formal groups}
\author{N.P. Strickland}
\email{N.P.Strickland@sheffield.ac.uk}

\address{Department of Pure Mathematics, University of 
Sheffield\\Sheffield S3 7RH, UK}



\begin{abstract}
 We show that a large class of formal groups can be realised
 functorially by even periodic ring spectra.  The main advance is in
 the construction of morphisms, not of objects.
\end{abstract}

\maketitle 

\section{Introduction}

Let $\FG$ be the category of formal groups (of the sort usually
considered in algebraic topology) over affine schemes.  Thus, an
object of $\FG$ consists of a pair $(G,S)$, where $S$ is an affine
scheme, $G$ is a formal group scheme over $S$, and a coordinate $x$
can be chosen such that $\OG\simeq\OS\psb{x}$ as $\OS$-algebras.  A
morphism from $(G_0,S_0)$ to $(G_1,S_1)$ is a commutative square
\[ \xymatrix {
    G_0 \rto^{\tp} \dto & G_1 \dto \\
    S_0 \rto_p & S_1
} \]
such that the induced map $G_0\xra{}p^*G_1$ is an isomorphism of
formal group schemes over $S_0$.

Next, recall that an \emph{even periodic ring spectrum} is a
commutative and associative ring spectrum $E$ such that $E^1=0$ and
$E^2$ contains a unit (which implies that $E\simeq\Sg^2E$ as spectra).
Here we are using the usual notation
$E^k=E^k(\text{point})=\pi_{-k}E$.  We write $\EPR$ for the category
of even periodic ring spectra.  (Everything here is interpreted in
Boardman's homotopy category of spectra; there are no $E_\infty$ or
$A_\infty$ structures.)

Given an even periodic ring spectrum $E$, we can form the scheme
$S_E:=\spec(E^0)$ and the formal group scheme $G_E=\spf(E^0\CPi)$ over
$S_E$.  This construction gives rise to a functor
$\Gm\:\EPR^\op\xra{}\FG$.

It is a natural problem to try to define a realisation functor
$R\:\FG\xra{}\EPR^\op$ with $\Gm R(G,S)\simeq(G,S)$, or at least to do
this for suitable subcategories of $\FG$.  For example, if we let
$\LFG$ denote the category of Landweber exact formal groups, and put
$\LEPR=\{E\in\EPR\st\Gm(E)\in\LFG\}$, one can show that the functor
$\Gm\:\LEPR^\op\xra{}\LFG$ is an equivalence; this is essentially due
to Landweber, but details of this formulation are given
in~\cite[Proposition 8.43]{st:fsfg}.  Inverting this gives a
realisation functor for $\LFG$, and many well-known spectra are
constructed using this.  In particular, this gives various different
versions of elliptic cohomology, based on various universal families
of elliptic curves over rings such as
$\Z[\sixth,c_4,c_6][\Delta^{-1}]$.

It is hard to say more than this unless we invert the prime $2$.  We
therefore make a blanket assumption:
\begin{assumption}
 From now on, all rings are assumed to be $\Z[\half]$-algebras.  In
 particular, we only consider schemes $S$ for which $2$ is invertible
 in $\OS$.  We use the symbol $MU$ for the spectrum that would
 normally be called $MU[\half]$.
\end{assumption}

The other main technique for constructing realisations is the
modernised version of Baas-Sullivan theory~\cite{ekmm:rma,st:pmm}.
This starts with a strictly commutative ring spectrum $R$, and an
algebra $A_*$ over $\pi_*R$, and it constructs a homotopically
commutative $R$-algebra spectrum $A$ with $\pi_*A=A_*$, provided that
$A_*$ has good structural properties.  Firstly, we assume as always
that $2$ is invertible in $A_*$.  Given this, the construction will
work if $A_*$ is a \emph{localised regular quotient (LRQ)} of $R_*$,
in other words it has the form $A_*=(S^{-1}\pi_*R)/I$, where $S$ is a
multiplicative set and $I$ is an ideal generated by a regular
sequence.  The construction can also be extended to cover the case
where $A_*$ is a free module over an LRQ of $\pi_*R$.

We can apply this taking $R$ to be the periodic bordism spectrum
\[ MP = \bigWedge_{n\in\Z}\Sg^{2n}MU[\half] \]
(we will verify in the appendix that this can be constructed as a
strictly commutative ring).  Given a formal group $(G,S)$ we can
choose a coordinate $x$, which gives a formal group law $F$ defined
over $\OS$, and thus a ring map $\pi_0MP\xra{}\OS$, making $\OS$ into
a $\pi_0MP$-algebra.  If this algebra has the right properties, then
we can use the Baas-Sullivan approach to construct $E$ with
$\Gm(E)\simeq(G,S)$.  It is convenient to make the following \emph{ad
  hoc} definition:
\begin{definition}
 A ring $R$ is \emph{standard} if $2$ is invertible in $R$ and $R$ is
 either a field or a ring of the form $T^{-1}\Z$ (for some set $T$ of
 primes).
\end{definition}
An easy argument given below shows that the above method can construct
realizations for all formal groups over standard rings.
Unfortunately, this construction is not obviously functorial: it
depends on a choice of coordinate, and morphisms of formal groups
do not generally preserve coordinates.  The main result of this paper
is to show that with suitable hypotheses we can nonetheless define a
functor. 

The basic point is to consider the situation where we have several
different coordinates, say $x_0,\ldots,x_r$ on a fixed formal group
$G$.  In a well-known way, this makes $\OS$ into an algebra over the
ring $\pi_0(MP^{(r+1)})$, and we can ask whether this can be realized
topologically by an $MP^{(r+1)}$-algebra; the question will be made
more precise in Section~\ref{sec-basic}.  We say that $G$ is
\emph{very good} if the question has an affirmative answer for all
$r\geq 0$ and all $x_0,\ldots,x_r$.
\begin{theorem}
 All formal groups over standard rings are very good.
\end{theorem}
This will be proved as Corollary~\ref{cor-very-good}.  

For our sharpest results, we need a slightly more complicated notion.
We say that a coordinate $x_0$ is \emph{multirealisable} if for any
list $x_1,\ldots,x_r$ of additional coordinates, the question
mentioned above has an affirmative answer.  We say that $G$ is
\emph{good} if it admits a multirealisable coordinate.  Of course, $G$
is very good iff \emph{every} coordinate is multirealisable.  We write
$\GFG$ for the category of good formal groups (considered as a
full subcategory of $\FG$).  The details are given in
Definition~\ref{defn-good}. 

\begin{theorem}\label{thm-LRQ-good}
 Let $x$ be a coordinate on a formal group $(G,S)$, and suppose that
 the classifying map $\pi_0MP\xra{}\OS$ makes $\OS$ into a localised
 regular quotient of $\pi_0MP$.  Then $x$ is multirealisable, and so
 $G$ is good.
\end{theorem}
This will be proved as Proposition~\ref{prop-LRQ-multi}.

\begin{corollary}
 At odd primes, the formal groups associated to $2$-periodic versions
 of $BP$, $P(n)$, $B(n)$, $E(n)$, $K(n)$, $k(n)$ and so on are all
 good.  \qed
\end{corollary}
This shows that there is a considerable overlap with the Landweber
exact case.  However, there are many good formal groups that are not
Landweber exact.  Conversely, there is no reason to expect that
Landweber exact formal groups will be good, although we have no
counterexamples. 

Our main result is as follows:
\begin{theorem}\label{thm-main}
 There is a realisation functor $R\:\GFG\xra{}\EPR$, with
 $\Gm R\simeq 1\:\GFG\xra{}\GFG$.
\end{theorem}
Note that good formal groups are realisable by definition; the content
of the theorem is that the realisation is well-defined and
functorial.

We next explain the formal part of the construction; in
Section~\ref{sec-proof} we will give additional details and prove that
we have the required properties.  The functor $R$ actually arises as
$UV^{-1}$ for a pair of functors $\GFG\xla{V}\CE\xra{U}\EPR$ in which
$V$ is an equivalence.  To explain $\CE$, recall that we have a
topological category $\Mod_0$ of $MP$-modules.  We write $\DMod_0$ for
the derived category, and $\EPA_0$ for the category of even periodic
commutative ring objects in $\DMod_0$.  The unit map $\eta\:S\xra{}MP$
gives a functor $\eta^*\:\EPA_0\xra{}\EPR$, and the objects of the
category $\CE$ are the objects $E\in\EPA_0$ for which the associated
coordinate on $\Gm(\eta^*E)$ is multirealisable.  The morphism set
$\CE(E_0,E_1)$ is a subset of $\EPR(\eta^*E_0,\eta^*E_1)$, the functor
$V\:\CE\xra{}\GFG$ is given by $\Gm$, and the functor
$U\:\CE\xra{}\EPR$ is given by $\eta^*$.  We say that a map
$f\:\eta^*E_0\xra{}\eta^*E_1$ in $\EPR$ is \emph{good} if there is a
commutative ring object $A$ in the derived category of $MP\Smash
MP$-modules together with maps $f'\:E_0\xra{}(\zin)^*A$ and
$f''\:(\zni)^*A\xra{}E_1$ in $\EPA_0$ such that $f''$ is an
equivalence and $f$ is equal to the composite
\[ \eta^*E_0 \xra{\eta^*f'} (\znn)^*A
    \xra{\eta^*f''} \eta^*E_1.
\]
The morphisms in the category $\CE$ are just the good maps.  To prove
Theorem~\ref{thm-main}, we need to show that
\begin{itemize}
 \item[(3)] The composite of two good maps is good, so $\CE$ really is
  a category.
 \item[(2)] For any map $\Gm(\eta^*E_0)\xra{}\Gm(\eta^*E_1)$ of good
  formal groups, there is a unique good map $\eta^*E_0\xra{}\eta^*E_1$
  inducing it, so that $V$ is full and faithful.  
 \item[(1)] For any good formal group $(G,S)$ there is an object
  $E\in\EPA_0$ such that $\Gm(\eta^*E)\simeq(G,S)$, so $V$ is
  essentially surjective.
\end{itemize}
To prove statement $(k)$, we need to construct modules over the
$k$-fold smash power of $MP$.  It will be most efficient to do this
for all $k$ simultaneously.

\section{Preliminaries}

\subsection{Differential forms}\label{subsec-forms}

Let $(G,S)$ be a formal group, and let $I\leq\OG$ be the augmentation
ideal.  Recall that the cotangent space of $G$ at zero is the module
$\om_G=I/I^2$.  If $x$ is a coordinate on $G$ that vanishes at zero,
then we write $dx$ for the image of $x$ in $I/I^2$, and note that
$\om_G$ is freely generated over $\OS$ by $dx$.  We define a graded
ring $D(G,S)^*$ by
\[ D(G,S)^k = \begin{cases}
               0                 & \text{ if $k$ is odd } \\
               \om_G^{\ot(-k/2)} & \text{ if $k$ is even. }
              \end{cases}
\]
Here the tensor products are taken over $\OS$, and $\om_G^{\ot n}$
means the dual of $\om_G^{\ot |n|}$ when $n<0$.  Where convenient, we
will convert to homological gradings by the usual rule:
$D(G,S)_k=D(G,S)^{-k}$. 

Now let $E$ be an even periodic ring spectrum with $\Gm(E)=(G,S)$.  We
then have $\OG=E^0\CPi$ and $I=\tE^0\CPi$ and one checks easily that
the inclusion $S^2=\CP^1\xra{}\CPi$ gives an isomorphism
$\om_G=I/I^2=\tE^0S^2=E^{-2}$.  Using the periodicity of $E$, we see
that this extends to a canonical isomorphism $D(\Gm(E))^*\simeq E^*$.

It also follows from this analysis (or from more direct arguments)
that a map $f\:E_0\xra{}E_1$ in $\EPR$ is a weak equivalence if and
only if $\pi_0f$ is an isomorphism.

\subsection{Periodic bordism}\label{subsec-MP}

Consider the homology theory $MP_*(X)=MU_*(X)\ot\Z[u,u^{-1}]$, where
$u$ has homological degree $2$ (and thus cohomological degree $-2$).
This is represented by the spectrum $MP=\bigWedge_{n\in\Z}\Sg^{2n}MU$,
with an evident ring structure.  It is well-known that $MU$ is an
$E_\infty$ ring spectrum; see for example~\cite[Section
IX]{lemast:esh}.  It is also shown there that $MU$ is an $H_\infty^2$
ring spectrum, which means (as explained in \cite[Remark
VII.2.9]{lemast:esh}) that $MP$ is an $H_\infty$ ring spectrum; this
is weaker than $E_\infty$ in theory, but usually equivalent in
practise.  As one would expect, $MP$ is actually an $E_\infty$ ring
spectrum; a proof is given in the appendix.  It follows
from~\cite[Proposition II.4.3]{ekmm:rma} that one can construct a
model for $MP$ that is a strictly commutative ring spectrum (or
``commutative $S$-algebra'').  We may also assume that it is a
cofibrant object in the category of all strictly commutative ring
spectra.

For typographical convenience, we write $MP(r)$ for the $(r+1)$-fold
smash power $MP\Smash\ldots\Smash MP$, which is again a strictly
commutative ring.  The spectra $MP(r)$ fit together into a
cosimplicial object in the usual way; for example, we have three maps
\[ \znii,\zini,\ziin \: MP(0) \xra{} MP(2). \]
In the category of strictly commutative ring spectra, the coproduct is
the smash product.  It follows formally that the smash product of
cofibrant objects is cofibrant, so in particular the objects $MP(r)$
are all cofibrant.

For $r>0$, it is well-known that $\pi_*MU^{(r+1)}$ is a polynomial
algebra over $\pi_*MU$ on countably many generators, and it follows
that there is a noncanonical isomorphism
\[ \pi_0MP(r) \simeq
    \pi_0MP[x_1,x_2,\ldots][x_1^{-1},\ldots,x_r^{-1}]. 
\]

There are $r+1$ obvious inclusions $MP\xra{}MP(r)$.  We can use these
to push forward the standard generator of $MP^0\CPi$, giving $r+1$
different coordinates on the formal group $\Gm(MP(r))$.  We denote
these by $\tx_0,\ldots,\tx_r$.

\subsection{Groups and laws}\label{subsec-laws}

We now define a category $\FG_r$ as follows.  The objects are systems
\[ (G,S,x_0,\ldots,x_r), \]
where $(G,S)$ is a formal group and the $x_i$ are coordinates on $G$.
The morphisms from $(G,S,x_0,\ldots,x_r)$ to $(H,T,y_0,\ldots,y_r)$
are the maps $(\tp,p)\:(G,S)\xra{}(H,T)$ in $\FG$ for which
$\tp^*y_i=x_i$ for all $i$.  Note that given $p$, the map $\tp$ is
determined by the fact that $\tp^*y_0=x_0$.  Thus, the forgetful
functor $(G,S,x_0,\ldots,x_r)\mapsto S$ (from $\FG_r$ to the category
of affine schemes) is faithful.

We also write $\Alg_r$ for the category of commutative algebras over
the ring $\pi_0MP(r)$.

\begin{proposition}\label{prop-FGr}
 There is an equivalence $\FG_r\simeq\Alg_r^{\op}$.
\end{proposition}
\begin{proof}
 Recall that we have coordinates $\tx_0,\ldots,\tx_r$ on
 $\Gm(MP(r))$.  Given an object $A\in\Alg_r$ we have a structure map
 $\spec(A)\xra{}\spec(\pi_0MP(r))$, and we can pull back $\Gm(MP(r))$
 to get a formal group $G_A$ over $\spec(A)$.  We can also pull back
 the coordinates $\tx_i$ to make $G_A$ an object of $\FG_r$.  It is
 easy to see that this construction defines a functor
 $U\:\Alg_r^{\op}\xra{}\FG_r$.  By forgetting down to the category of
 affine schemes, we see that $U$ is faithful.

 We now claim that $U$ is an equivalence.  We will deduce this from a
 well-known result of Quillen by a sequence of translations.  First,
 Quillen tells us that maps $\pi_*MU^{(r+1)}\xra{}B_*$ of graded rings
 biject naturally with systems 
 \[ F_0 \xla{f_0} F_1 \xla{f_1} \cdots \xla{f_{r-1}} F_r, \]
 where each $F_i$ is a homogeneous formal group law over $B_*$ and
 each $f_i$ is a strict isomorphism.  By a standard translation to the
 even periodic case, we see that maps $\pi_0MP(r)\xra{}A$ of ungraded
 rings biject naturally with systems 
 \[ F_0 \xla{f_0} F_1 \xla{f_1} \cdots \xla{f_{r-1}} F_r, \]
 where each $F_i$ is a formal group law over $A$ and each $f_i$ is a
 (not necessarily strict) isomorphism.

 Now suppose we have an object $(G,S,x_0,\ldots,x_r)$ in $\FG_r$.  For
 each $i$ there is a unique formal group law $F_i$ over $\OS$ such
 that $x_i(a+b)=F_i(x_i(a),x_i(b))$ for sections $a,b$ of $G$.
 Moreover, as $x_{i+1}$ is another coordinate, we can write
 $x_i=f_i(x_{i+1})$ for a unique power series $f_i\in\OS\psb{t}$.  It
 is easy to check that $f_i$ is an isomorphism from $F_{i+1}$ to
 $F_i$, so Quillen's theorem gives us a map $\pi_0MP(r)\xra{}\OS$,
 allowing us to regard $\OS$ as an object of $\Alg_r$.  It is easy to
 see that this construction gives a functor
 $\FG_r\xra{}\Alg_r^{\op}$.  We leave it to the reader to check that
 this is inverse to $U$.
\end{proof}

\subsection{Module categories}

We write $\Mod_r$ for the category of $MP(r)$-modules (in the strict
sense, not the homotopical one).  Note that a map $f\:A_0\xra{}A_1$ of
strictly commutative ring spectra gives a functor
$f^*\:\Mod_{A_1}\xra{}\Mod_{A_0}$, which is just the identity on the
underlying spectra (and thus preserves weak equivalences).  It follows
easily that for any two maps $A_0\xra{f}A_1\xra{g}A_2$, the functor
$f^*g^*$ is actually equal (not just naturally isomorphic or naturally
homotopy equivalent) to $(gf)^*$.  Thus, the categories $\Mod_r$ fit
together to give a simplicial category $\Mod_*$.

\begin{remark}
 For us, a \emph{simplicial category} means a simplicial object in the
 category of categories.  Elsewhere in the literature, the same phrase
 is sometimes used to refer to categories enriched over the category
 of simplicial sets, which is a rather different notion.
\end{remark}

Next, we write $\DMod_r$ the derived category of $\Mod_r$, as
in~\cite[Chapter III]{ekmm:rma}.  As usual, there are two different
models for a category such as $\DMod_r$:
\begin{itemize}
 \item[(a)] One can take the objects to be the cofibrant objects in
  $\Mod_r$, and morphisms to be homotopy classes of maps; or
 \item[(b)] One can use all objects in $\Mod_r$ and take morphisms to
  be equivalence classes of ``formal fractions'', in which one is
  allowed to invert weak equivalences.
\end{itemize}
We will use model~(b).  This preserves the strong functorality
mentioned previously, and ensures that $\DMod_*$ is again a
simplicial category.

We also write $\EPA_r$ for the category of even periodic commutative
ring objects in $\DMod_r$, giving another simplicial category.  (Note
that periodicity is actually automatic, because $MP(r)$ is itself
periodic.)  Various fragments of the simplicial structure will be used
in Section~\ref{sec-proof}.

\section{Basic realisation results}
\label{sec-basic}

Let $R$ be a strictly commutative ring spectrum that is even and
periodic, such that $R_0$ is an integral domain (and as always, $2$ is
invertible).  The main examples will be $R=MP(r)$ for $r\geq 0$.  Let
$\mathcal{D}$ be the derived category of $R$-modules, and let $\CR$ be
the category of commutative ring objects $A\in\mathcal{D}$ such that
$\pi_1A=0$.  Recall that if $f$ is a morphism in $\CR$ such that
$\pi_0f$ is an isomorphism, then $\pi_*f$ is also an isomorphism and
so $f$ is an equivalence.

We also write $\CR_0$ for the category of commutative
algebras over $\pi_0R$.  We say that an object $A\in\CR$ is
\emph{strong} if for all $B\in\CR$, the map 
\[ \pi_0\: \CR(A,B) \xra{} \CR_0(\pi_0A,\pi_0B) \]
is a bijection.  A \emph{realisation} of an object $A_0\in\CR_0$ is a
pair $(A,u)$, where $A\in\CR$ and $u\:\pi_0A\xra{}A_0$ is an
isomorphism.  We say that $(A,u)$ is a \emph{strong realisation} iff
the object $A$ is strong; if so, we have a natural isomorphism
$\CR(A,B)\simeq\CR_0(A_0,\pi_0B)$.  We say that $A_0$ is
\emph{strongly realisable} if it admits a strong realisation.  If so,
it is easy to check that all realisations are strong, and any two
realisations are linked by a unique isomorphism.

The results of~\cite{st:pmm} provide a good supply of strongly
realisable algebras, except that we need a little translation between
the even periodic framework and the usual graded framework.  Suppose
that $A_0\in\CR_0$, and put $T=\spec(A_0)$.  We have a unit map
$\eta\:\pi_0R\xra{}A_0$ and thus a map $\spec(\eta)\:T\xra{}S_R$; we
can pull back the formal group $G_R$ along this to get a formal group
$H:=\spec(\eta)^*G_R$ over $T$.  From this we get a map
$\eta_*\:R_*=D(G_R,S_R)_*\xra{}D(H,T)_*$, which agrees with $\eta$ in
degree zero.  Indeed, if we choose a generator $u$ of $R_2$ over
$R_0$, then $\eta_*$ is just the map
$R_0[u,u^{-1}]\xra{}A_0[u,u^{-1}]$ obtained in the obvious way from
$\eta$.  It is easy to check that $A_0$ is strongly realisable (as
defined in the previous paragraph) iff $D(H,T)_*$ is strongly
realisable over $R_*$ (as defined in~\cite{st:pmm}).

\begin{definition}
 A \emph{short ordinal} is an ordinal $\lm$ of the form $n.\om+m$ for
 some $n,m\in\N$.  A \emph{regular sequence} in a ring $R_0$ is a
 system of elements $(x_\al)_{\al<\lm}$ for some short ordinal $\lm$
 such that $x_\al$ is not a zero-divisor in the ring
 $(S^{-1}R_0)/(x_\bt\st\bt<\al)$.  An object $A_0\in\CR_0$ is a
 \emph{localised regular quotient} (or LRQ) of $R_0$ if
 $A_0=(S^{-1}R_0)/I$ for some subset $S\subset R_0$ and some ideal
 $I\leq S^{-1}R_0$ that can be generated by a regular sequence.
\end{definition}
\begin{remark}
 We have made a small extension of the usual notion of a regular
 sequence, to ensure that any LRQ of an LRQ of $R_0$ is itself an LRQ
 of $R_0$; see Lemma~\ref{lem-LRQ-LRQ}.
\end{remark}

\begin{proposition}\label{thm-LRQ}
 If $A_0$ is an LRQ of $R_0$, then it is strongly realisable.
\end{proposition}
\begin{proof}
 This is essentially~\cite[Theorem 2.6]{st:pmm}, translated into a
 periodic setting as explained above.  Here we are using a slightly
 more general notion of a regular sequence, but all the arguments can
 be adapted in a straightforward way.  The main point is that any
 countable limit ordinal has a cofinal sequence, so homotopy colimits
 can be constructed using telescopes in the usual way.  Andrey Lazarev
 has pointed out a lacuna in~\cite{st:pmm}: it is necessary to assume
 that the elements $x_\al$ are all regular in $S^{-1}R_0$ itself,
 which is not generally automatic.  However, we are assuming that
 $R_0$ is an integral domain so this issue does not arise.
\end{proof}

\begin{proposition}\label{prop-tensor}
 Suppose that
 \begin{itemize}
  \item $A$ and $B$ are strong realisations of $A_0$ and $B_0$
  \item The natural map $A_0\ot_{R_0}B_0\xra{}(A\Smash_R B)_0$ is an
   isomorphism. 
 \end{itemize}
 Then $A\Smash_RB$ is a strong realisation of $A_0\ot_{R_0}B_0$.
\end{proposition}
\begin{proof}
 This follows from~\cite[Corollary 4.5]{st:pmm}.
\end{proof}

\begin{proposition}\label{prop-free}
 If $A_0\in\CR_0$ is strongly realisable, and $B_0$ is an algebra
 over $A_0$ that is free as a module over $A_0$, then $B_0$ is also
 strongly realisable.
\end{proposition}
\begin{proof}
 This follows from~\cite[Proposition 4.13]{st:pmm}.
\end{proof}

\begin{proposition}\label{prop-Z-ninv}
 Suppose that $R_0$ is a polynomial ring in countably many variables
 over $\Z[\half]$, that $A_0\in\CR_0$, and that $A_0=\Z[1/2n]$ as a
 ring (for some $n$).  Then $A_0$ is an LRQ of $R_0$, and thus is
 strongly realisable.
\end{proposition}
\begin{proof}
 Choose a system of polynomial generators $\{x_k\st k\geq 0\}$ for
 $R_0$ over $\Z[\half]$.  Put $a_k=\eta(x_k)\in A_0=\Z[1/n]$ and
 $y_k=x_k-a_k\in R_0[1/n]$.  It is clear that
 $R_0[1/2n]=\Z[1/2n][y_k\st k\geq 0]$, that the elements $y_k$ form a
 regular sequence generating an ideal $I$ say, and that
 $A_0=R_0[1/2n]/I$.
\end{proof}

\begin{proposition}\label{prop-field}
 Suppose that $R_0$ is a polynomial ring in countably many variables
 over $\Z[\half]$, that $A_0\in\CR_0$, and that $A_0$ is a field
 (necessarily of characteristic different from $2$).  Then $A_0$ is a
 free module over an LRQ of $R_0$, and thus is strongly realisable.
\end{proposition}
\begin{proof}
 For notational simplicity, we assume that $A_0$ has characteristic
 $p>2$; the case of characteristic $0$ is essentially the same.

 Choose a set $X$ of polynomial generators for $R_0$ over $\Z[\half]$.
 Let $K$ be the subfield of $A_0$ generated by the image of $\eta$, or
 equivalently by $\eta(X)$.  We can choose a subset $Y\sse X$ such
 that $\eta(Y)$ is a transcendence basis for $K$ over $\Fp$.  This
 means that the subfield $L_0$ of $K$ generated by $\eta(Y)$ is
 isomorphic to the rational function field $\Fp(Y)$, and that $K$ is
 algebraic over $L_0$.  Put $S=\Z[\half,Y]\sm (p\Z[\half,Y])$, so
 $L_0=(S^{-1}\Z[\half,Y])/p$.  Next, list the elements of $X\sm Y$ as
 $\{x_1,x_2,\ldots\}$, and let $L_k$ be the subfield of $K$ generated
 by $\{x_i\st i\leq k\}$.  (We will assume that $X\sm Y$ is infinite;
 if not, the notation changes slightly.)  As $x_k$ is algebraic over
 $L_{k-1}$, there is a monic polynomial $f_k(t)\in L_{k-1}[t]$ with
 $L_k=L_{k-1}[x_k]/f_k(x_k)$.  As $L_{k-1}$ is a quotient of the ring
 $P_{k-1}:=S^{-1}\Z[Y,x_1,\ldots,x_{k-1}]$, we can choose a monic
 polynomial $g_k(t)\in P_{k-1}[t]$ lifting $f_k$, and put
 $z_k:=g_k(x_k)\in P_k\sse S^{-1}R_0$.  It is not hard to check that
 the sequence $(p,z_1,z_2,\ldots)$ is regular in $S^{-1}R_0$, and that
 $(S^{-1}R_0)/(z_i\st i>0)=K$, so $K$ is an LRQ of $R_0$.  It is clear
 that $A_0$ is free over the subfield $K$.
\end{proof}

\begin{lemma}\label{lem-LRQ-LRQ}
 An LRQ of an LRQ is an LRQ.
\end{lemma}
\begin{proof}
 Suppose that $B=(S^{-1}A)/(x_\al\st\al<\lm)$ and
 $C=(T^{-1}B)/(y_\bt\st\bt<\mu)$, where $\lm$ and $\mu$ are short
 ordinals, and the $x$ and $y$ sequences are regular in $S^{-1}A$ and
 $T^{-1}B$ respectively.  Let $T'$ be the set of elements of $A$ that
 become invertible in $T^{-1}B$; clearly $S\sse T'$ and
 $T^{-1}B=((T')^{-1}A)/(x_\al\st\al<\lm)$.  As $(T')^{-1}A$ is a
 localisation of $S^{-1}A$ and localisation is exact, we see that $x$
 is a regular sequence in $(T')^{-1}A$ as well.  After multiplying by
 suitable elements of $T'$ if necessary, we may assume that $y_\bt$
 lies in the image of $A$ (this does not affect regularity, as the
 elements of $T'$ are invertible).  We then put $z_\al=x_\al$ for
 $\al<\lm$, and let $z_{\lm+\bt}$ be any preimage of $y_\bt$ in $A$
 for $0\leq\bt<\mu$.  This gives a regular sequence in $(T')^{-1}A$
 indexed by $\lm+\mu$, such that
 $C=((T')^{-1}A)/(z_\gm\st\gm<\lm+\mu)$ as required.
\end{proof}

We now specialize to the case $R=MP(r)$, so $\CR_0=\EPA_r$.  We write
$\Gm_r$ for the evident composite functor
\[ \EPA_r^{\op} \xra{\pi_0} \Alg_r^{\op} \simeq \FG_r. \]
Translating our previous definitions via the equivalence
$\Alg_r^{\op}\simeq\FG_r$, we obtain the following.
\begin{definition}\label{defn-strong-FGr}
 An object $A\in\EPA_r$ is \emph{strong} if for all $B\in\EPA_r$, the
 map
 \[ \Gm_r \: \EPA_r(A,B) \xra{} \FG_r(\Gm_r(B),\Gm_r(A)) \]
 is a bijection.  
\end{definition}
\begin{definition}
 A \emph{realisation} of an object $(G,S,\ux)\in\FG_r$ is a triple
 $(A,\tp,p)$, where $A\in\EPA_r$ and $(\tp,p)\:\Gm_rA\xra{}(G,S,\ux)$
 is an isomorphism.  This is a \emph{strong realisation} if the object
 $A$ is strong.
\end{definition}

We now give more precise versions of the definitions in the
introduction. 
\begin{definition}
 A formal group $(G,S)$ is \emph{very good} if for every nonempty list
 $\ux$ of coordinates, the object $(G,S,\ux)\in\FG_r$ is strongly
 realisable. 
\end{definition}

\begin{definition}\label{defn-good}
 A coordinate $x_0$ on $G$ is \emph{multirealisable} if for every list
 $x_1,\ldots,x_r$ of coordinates, the object
 $(G,S,x_0,\ldots,x_r)\in\FG_r$ is strongly realisable.  A formal
 group $(G,S)$ is \emph{good} if it admits a multirealisable
 coordinate.  We write $\GFG$ for the category of good formal groups.
\end{definition}
\begin{remark}\label{rem-perm}
 Let $x_0,\ldots,x_r$ be coordinates, and suppose that $x_0$ is
 multirealisable.  Let $\sg$ be a permutation of $\{0,\ldots,r\}$.
 Using the evident action of permutations on $MP(r)$, we see that the
 object $(G,S,x_{\sg(0)},\ldots,x_{\sg(r)})$ is strongly realisable.
\end{remark}

\begin{proposition}\label{prop-LRQ-multi}
 Suppose that $x_0$ is such that the classifying map
 $\pi_0MP\xra{}\OS$ makes $\OS$ an LRQ of $\pi_0MP$.  Then $x_0$ is
 multirealisable, so $(G,S)$ is good.
\end{proposition}
\begin{proof}
 The coordinate $x_0$ gives a map $f_0\:\pi_0MP\xra{}\OS$.  By
 assumption, there is a multiplicative set $T\sse\pi_0MP$ and a
 regular ideal $I$ such that $f_0$ induces an isomorphism
 $(T^{-1}\pi_0MP)/I\xra{}\OS$.  

 Now consider a list of additional coordinates $x_1,\ldots,x_r$ say.
 These give a map $f\:\pi_0MP(r)\xra{}\OS$ extending $f_0$.  We know
 from Section~\ref{subsec-MP} that $\pi_0MP(r)$ is a polynomial ring
 in countably many variables over $\pi_0MP$, in which $r$ of the
 variables have been inverted, so we can write
 \[ \pi_0MP(r) = \pi_0MP[u_1,u_2,\ldots][u_1^{-1},\ldots,u_r^{-1}]. \]
 Put 
 \[ A_0 = \OS[u_1,u_2,\ldots][u_1^{-1},\ldots,u_r^{-1}], \]
 which is evidently an LRQ of $\pi_0MP(r)$.  It is easy to see that
 $f$ induces a map $f'\:A_0\xra{}\OS$ of $\OS$-algebras.  Put
 $a_k=f'(u_k)\in\OS$, and $v_k=u_k-a_k\in A_0$.  Clearly $A_0$ is a
 localisation of $\OS[v_k\st k>0]$, the sequence of $v$'s is regular
 in $A_0$, and $A_0/(v_k\st k>0)=\OS$ as $\pi_0MP(r)$-algebras.  It
 follows that $\OS$ is an LRQ of an LRQ, and thus an LRQ, over
 $\pi_0MP(r)$.  It is thus strongly realisable as required.
\end{proof}
\begin{corollary}\label{cor-very-good}
 If $\OS$ is a standard ring, then every coordinate is
 multirealisable, and so $(G,S)$ is very good.
\end{corollary}
\begin{proof}
 This now follows from Propositions~\ref{prop-Z-ninv}
 and~\ref{prop-field}. 
\end{proof}

\section{Proof of the main theorem}
\label{sec-proof}

Let $\CE$ denote the class of objects $E\in\EPA_0$ for which the
resulting coordinate on $\Gm(\eta^*E)$ is multirealisable.  Note that
this means that $\Gm_1E$ is strongly realisable, so every realisation
is strong, so in particular $E$ is a strong object.

\begin{proposition}\label{prop-surj}
 For any good formal group $(G,S)$, there exists $E\in\CE$ with
 $\Gm(\eta^*E)\simeq(G,S)$.
\end{proposition}
\begin{proof}
 By the definition of goodness we can choose a multirealisable
 coordinate $x_0$ on $G$.  This means in particular that the object
 $(G,S,x_0)\in\FG_0$ is isomorphic to $\Gm_0(E)$ for some
 $E\in\EPA_0$.  It follows that $(G,S)\simeq\Gm(\eta^*E)$, as
 required.  
\end{proof}

\begin{proposition}\label{prop-maps}
 Suppose we have objects $E_0,E_1\in\CE$, together with a map
 \[ (\tp,p)\:\Gm(\eta^*E_1)\xra{}\Gm(\eta^*E_0) \]
 in $\GFG$.  Then there is a unique good map
 $f\:\eta^*E_0\xra{}\eta^*E_1$ such that $\Gm(f)=(\tp,p)$.
\end{proposition}
\begin{proof}
 We first put $(G_i,S_i,x_i)=\Gm_0E_i$ for $i=0,1$.  

 We introduce a category $\CB=\CB(E_0,E_1,\tp,p)$ as follows.  The
 objects are triples $(A,f',f'')$ where
 \begin{itemize}
  \item[(a)] $A$ is an object of $\EPA_1$.
  \item[(b)] $f'$ is a morphism $E_0\xra{}(\zin)^*A$ in
   $\EPA_0$. 
  \item[(c)] $f''$ is an isomorphism $(\zni)^*A\xra{}E_1$ in
   $\EPA_0$. 
  \item[(d)] The composite
   \[ f = \tht(A,f',f'') := 
       (\eta^*E_0\xra{\eta^*f'}(\znn)^*A\xra{\eta^*f''}\eta^*E_1)
   \]
   satisfies $\Gm(f)=(\tp,p)$.
 \end{itemize}
 The morphisms from $(A,f',f'')$ to $(B,g',g'')$ in $\CB$ are the
 isomorphisms $u\:A\xra{}B$ in $\EPA_1$ for which
 $((\zin)^*u)f'=g'$ and $g''((\zni)^*u)=f''$.  

 The maps of the form $\tht(A,f',f'')$ are precisely the good maps
 that induce $(\tp,p)$, and isomorphic objects of $\CB$ have the same
 image under $\tht$.  It will thus suffice to show that
 $\CB\neq\emptyset$ and all objects of $\CB$ are isomorphic.

 First, as $x_1$ is multirealisable, we can find an object
 $A\in\EPA_1$ and an isomorphism
 $(\tq,q)\:\Gm_1A\xra{}(G_1,S_1,\tp^*x_0,x_1)$ displaying $A$ as a
 strong realisation of $(G_1,S_1,\tp^*x_0,x_1)$.  We write
 $(H,T,y_0,y_1)=\Gm_1A$, so $(\tq,q)\:(H,T)\xra{\simeq}(G_1,S_1)$ and
 $(\tp\tq)^*x_0=y_0$ and $\tq^*x_1=y_1$.  We can thus regard
 $(\tp\tq,pq)$ as a morphism
 \[ \Gm_0((\zin)^*A) = (H,T,y_0) \xra{} (G_0,S_0,x_0) = \Gm_0E_0, \]
 and $E_0$ is a strong realisation of $(G_0,S_0,x_0)$, so this must
 come from a map $f'\:E_0\xra{}(\zin)^*A$ in $\EPA_0$.  Similarly, we
 can regard $(\tq,q)$ as an isomorphism 
 \[ \Gm_0((\zni)^*A) = (H,T,y_1) \xra{\simeq}
     (G_1,S_1,x_1) = \Gm_0^*E_1.
 \]
 As $E_1$ is a strong realisation of $(G_1,S_1,x_1)$, this comes from
 a map $E_1\xra{}(\zni)^*A$; this is easily seen to be an isomorphism,
 and we let $f''\:(\zni)^*A\xra{}E_1$ be the inverse map.  It is then
 clear that the map
 \[ f=(\eta^*f'')\circ(\eta^*f')\:\eta^*E_0\xra{}\eta^*E_1 \] 
 is good and satisfies $\Gm(f)=(\tp,p)$, so $(A,f',f'')\in\CB$.  Thus
 $\CB\neq\emptyset$.  

 Now suppose we have another object $(B,g',g'')\in\CB$, with
 $\Gm_1B=(K,U,z_0,z_1)$ say.  We put
 \begin{align*}
  (\tr,r) = \Gm_1g'' &\: (G_1,S_1,x_1) \xra{\simeq} 
     \Gm_1((\zni)^*B) = (K,U,z_1)  \\
  (\ts,s) = \Gm_1g'  &\: \Gm_1((\zin)^*B) = (K,U,z_0)
      \xra{} (G_0,S_0,x_0).
 \end{align*}
 By hypothesis we have
 $(\ts\tr,sr)=(\tp,p)\:(G_1,s_1)\xra{}(G_0,S_0)$.  We display all
 these maps in the following commutative diagram:
 \[ \xymatrix @=3pc {
   (H,T) \rto^{(\tq,q)}_{\simeq} \dto_{(\tp\tq,pq)} &
   (G_1,S_1) \dlto_{(\tp,p)} \dto_{\simeq}^{(\tr,r)} \\
   (G_0,S_0) &
   (K,U) \lto^{(\ts,s)}.
 } \]
 We now claim that $(\tr\tq,rq)$ can be regarded as a map
 \[ (H,T,y_0,y_1)\xra{}(K,U,z_0,z_1). \]
 Indeed, it is clear from the data recorded above that it is a map
 $(H,T,y_1)\xra{}(K,U,z_1)$, so it will suffice to check that
 $(\tr\tq)^*z_0=y_0$.  We are given that $z_0=\ts^*x_0$ and
 $\ts\tr=\tp$ and $(\tp\tq)^*x_0=y_0$; the claim follows.  As $r$ and
 $q$ are isomorphisms, we have an isomorphism
 \[ (\tr\tq,rq)^{-1}\: \Gm_1B = (K,U,z_0,z_1) \xra{}
                       (H,T,y_0,y_1) = \Gm_1A
 \]
 in $\FG_1$.  As $A$ is a strong realization, this comes from a unique
 map $u\:A\xra{}B$ in $\EPA_1$, which is easily seen to be an
 isomorphism. 

 We must show that $u$ is a morphism in our category $\CB$, or
 equivalently that in $\EPA_0$ we have
 \begin{align*}
  ((\zin)^*u)f'=g'   & \: E_0\xra{}(\zin)^*B \\
  g''((\zni)^*u)=f'' & \: (\zni)^*B\xra{} E_1.  
 \end{align*}
 Note that $E_0$ and $E_1$ are strong, and $f''$ is an isomorphism, so
 $(\zni)^*B$ is strong.  It is thus enough to check our two equations
 after applying $\pi_0$ (here we have used the original definition
 rather than the equivalent one in Definition~\ref{defn-strong-FGr}).  
 By definition or construction, we have
 \begin{align*}
  \spec(\pi_0f')  &= pq \\
  \spec(\pi_0f'') &= q^{-1} \\
  \spec(\pi_0g')  &= s \\
  \spec(\pi_0g'') &= r \\
  \spec(\pi_0u)   &= (rq)^{-1} \\
  sr              &= p.
 \end{align*}
 It follows easily that $(\pi_0u)(\pi_0f')=\pi_0g'$ and
 $(\pi_0g'')(\pi_0u)=\pi_0f''$, as required.  This shows that $u$ is
 an isomorphism in $\CB$, and thus that $f$ is the unique good map
 inducing the map $(\tp,p)$.
\end{proof}

\begin{lemma}\label{lem-id}
 For any $E\in\CE$, the identity map $1\:\eta^*E\xra{}\eta^*E$ is
 good. 
\end{lemma}
\begin{proof}
 Note that the multiplication map $MP(1)=MP\Smash MP\xra{}MP$ is a map
 of ring spectra (in the strict sense) and so induces a functor
 $\mu^*\:\EPA_0\xra{}\EPA_1$ with $(\zin)^*\mu^*E=(\zni)^*\mu^*E=E$ on
 the nose.  We can thus take $A=\mu^*E$ and $f'=f''=1_E$ to show that
 $1_E$ is good.
\end{proof}

\begin{proposition}\label{prop-comp}
 Suppose we have objects $E_0,E_1,E_2\in\CE$ and good morphisms
 $\eta^*E_0\xra{f}\eta^*E_1\xra{g}\eta^*E_2$.  Then the composite $gf$
 is also good.
\end{proposition}
\begin{proof}
 Write $(G_i,S_i,x_i)=\Gm_0E_i$ and
 $(\tp,p)=\Gm(f)\:(G_1,S_1)\xra{}(G_0,S_0)$ and
 $(\tq,q)=\Gm(g)\:(G_2,S_2)\xra{}(G_1,S_1)$.

 Choose objects $A,B\in\EPA_1$ and maps
 \begin{align*}
  f'  &\: E_0 \xra{} (\zin)^*A \\
  f'' &\: (\zni)^*A \xra{\simeq} E_1 \\
  g'  &\: E_1 \xra{} (\zin)^*B \\
  g'' &\: (\zni)^*B \xra{\simeq} E_2
 \end{align*}
 exhibiting the goodness of $f$ and $g$.  This gives rise to
 isomorphisms 
 \begin{align*}
  \Gm_1A &= (G_1,S_1,\tp^*x_0,x_1) \\
  \Gm_1B &= (G_2,S_2,\tq^*x_1,x_2).
 \end{align*}

 Next, observe that we have an object
 $(G_2,S_2,(\tp\tq)^*x_0,\tq^*x_1,x_2)\in\FG_2$, which is strongly
 realisable because $x_2$ is a multirealisable coordinate.  We can
 thus choose an object $P\in\EPA_2$ and an isomorphism
 \[ (\tr,r)\: \Gm_2P\xra{}(G_2,S_2,(\tp\tq)^*x_0,\tq^*x_1,x_2) \]
 making $P$ a strong realisation.  We can also regard $(\tr,r)$ as an
 isomorphism 
 \[ \Gm_1((\znii)^*P) \xra{} (G_2,S_2,\tq^*x_1,x_2) = \Gm_1B. \]
 As $B$ is strong, this comes from a unique isomorphism
 $v\:(\znii)^*P\xra{}B$ in $\EPA_1$.

 Similarly, we can regard $(\tr,r)$ as an isomorphism
 \[ \Gm_1((\ziin)^*P)\xra{}(G_2,S_2,\tq^*\tp^*x_0,\tq^*x_1), \]
 and we can regard $(\tq,q)$ as a morphism
 \[ (G_2,S_2,\tq^*\tp^*x_0,\tq^*x_1) \xra{}
     (G_1,S_1,\tp^*x_0,x_1) \simeq \Gm_1A.
 \]
 As $A$ is strong, the composite $(\tq\tr,qr)$ must come from a unique
 map $u\:A\xra{}(\ziin)^*P$ in $\EPA_1$.
 
 We now put 
 \begin{align*}
  C  &= (\zini)^*P\in\EPA_1 \\
  h' &= (E_0 \xra{f'} (\zin)^*A\xra{(\zin)^*u} (\zinn)^*P= (\zin)^*C) \\
  h''&= ((\zni)^*C = (\znni)^*P\xra{(\zni)^*v}(\zni)^*B \xra{g''} E_2).
 \end{align*}
 As $v$ and $g''$ are isomorphisms, the same is true of $h''$.  We
 claim that after forgetting down to $\EPR$, we have $h''h'=gf$; this
 will prove that $gf$ is good as claimed.  We certainly have
 $h''h'=g''vuf'$ and $gf=g''g'f''f'$ so it will suffice to show that
 $vu=g'f''\:A\xra{}B$ in $\EPR$.  For this, it will be enough to prove
 that the following diagram in $\EPA_0$ commutes.
 \[ \xymatrix {
   (\zni)^*A \rto^{(\zni)^*u} \dto_{f''}^{\simeq} & (\znin)^*P
   \dto_{\simeq}^{(\zin)^*v} \\ E_1 \rto_{g'} & (\zin)^*B. 
 } \]
 As this is a diagram in $\EPA_0$ and $(\zni)^*A\simeq E_1$ is strong,
 it will be enough to check that the diagram commutes after applying
 $\pi_0$.  By construction we have
 $\pi_0(u)=w^{-1}\circ\psi\circ\pi_0(f'')$ and
 $\psi=\pi_0(g)=\pi_0(g'')\circ\pi_0(g')$ and
 $\pi_0(v)=\pi_0(g'')^{-1}\circ w$.  It follows directly that the above
 diagram commutes on homotopy, groups, so it commutes in $\EPA_0$, so
 it commutes in $\EPR$, so $gf=h''h'$ in $\EPR$ as explained
 previously.  Thus, the map $gf$ is good, as claimed.
\end{proof}

\begin{proof}[Proof of Theorem~\ref{thm-main}]
 We merely need to collect results together and explain the argument
 in the introduction in more detail.  Lemma~\ref{lem-id} and
 Proposition~\ref{prop-comp} show that we can make $\CE$ into a
 category by taking the good maps from $\eta^*E_0$ to $\eta^*E_1$ as
 the morphisms from $E_0$ to $E_1$.  Tautologically, we can define a
 faithful functor $U\:\CE\xra{}\EPR$ by $U(E)=\eta^*E$ and $U(f)=f$.
 We then define $V=\Gm U\:\CE\xra{}\FG$; by the definition of $\CE$,
 this actually lands in $\GFG$.  Proposition~\ref{prop-surj} says that
 $V$ is essentially surjective, and Proposition~\ref{prop-maps} says
 that $V$ is full and faithful.  This means that $V$ is an
 equivalence, so we can invert it and define
 $R=UV^{-1}\:\GFG\xra{}\EPR$.  As $V=\Gm U$ we have $\Gm R=1$, so $R$
 is the required realisation functor.
\end{proof}

\appendix
\section{Appendix : The product on $MP$}

In this appendix we verify that $MP$ can be constructed as an
$E_\infty$ ring spectrum.

Let $\CU$ be a complex universe.  For any finite-dimensional subspace
$U$ of $\CU$, we write $U_L=U\oplus 0<\CU\oplus\CU$ and 
$U_R=0\oplus U<\CU\oplus\CU$.  We let $\Grass(U\oplus U)$ denote the
Grassmannian of all subspaces of $U\oplus U$ (of all possible
dimensions), and we write $\gm_U$ for the tautological bundle over
this space, and $\Thom(U\oplus U)$ for the associated Thom space.  If
$U\leq U'<\CU$ then we define $i\:\Grass(U^2)\xra{}\Grass((U')^2)$
by $i(A)=A\oplus(U'\ominus U)_R$.  On passing to Thom spaces we get a map
$\sg\:\Sg^{U'\ominus U}\Thom(U^2)\xra{}\Thom((U')^2)$.  These maps
can be used to assemble the spaces $\Thom(U^2)$ into a
$\Sg$-inclusion prespectrum indexed by the complex subspaces of $\CU$.
We write $T_\CU$ for this prespectrum, and $MP_\CU$ for its
spectrification.

Now let $\CV$ be another complex universe, so we have a prespectrum
$T_\CV$ over $\CV$, and thus an external smash product 
$T_\CU\Smash_\ext T_\CV$ indexed on the complex subspaces of
$\CU\oplus\CV$ of the form $U\oplus V$.  The direct sum gives a map
$\Grass(U^2)\tm\Grass(V^2)\xra{}\Grass((U\oplus V)^2)$
which induces a map
$\Thom(U^2)\Smash\Thom(V^2)\xra{}\Thom((U\oplus V)^2)$.
These maps fit together to give a map 
$T_\CU\Smash_\ext T_\CV\xra{}T_{\CU\oplus\CV}$, and thus a map
$MP_\CU\Smash_\ext MP_\CV\xra{}MP_{\CU\oplus\CV}$ of spectra over
$\CU\oplus\CV$.  Essentially the same construction gives maps 
\[ MP_{\CU_1}\Smash_\ext\ldots\Smash_\ext MP_{\CU_r}
    \xra{} MP_{\CU_1\oplus\ldots\oplus\CU_r}.
\]
If $\CU_1=\ldots\CU_r=\CU$, then this map is $\Sg_r$-equivariant.

Now suppose instead that we have a complex linear isometry
$f\:\CU\xra{}\CV$.  This gives evident homeomorphisms
$\Thom(U^2)\xra{}\Thom((fU)^2)$, which fit together to induce a map
$MP_\CU\xra{}f^*MP_\CV$, which is adjoint to a map
$f_*MP_\CU\xra{}MP_\CV$.  We next observe that this construction is
continuous in all possible variables, including $f$.  (This statement
requires some interpretation, but there are no new issues beyond those
that are well-understood for $MU$; the cleanest technical framework is
provided by~\cite{el:ggs}.)  It follows that they fit together to give
a map $\CLC(\CU,\CV)\ths MP_\CU\xra{}MP_\CV$ of spectra over $\CV$.

We now combine this with the product structure mentioned earlier to
get a map
\[ \CLC(\CU^r,\CU)\ths_{\Sg_r}
    (MP_\CU\Smash_\ext\ldots\Smash_\ext MP_\CU) \xra{}
     MP_\CU.
\]
This means that $MP_\CU$ has an action of the $E_\infty$ operad of
complex linear isometries, as required.

All that is left is to check that the spectrum $MP=MP_{\C^\infty}$
constructed above has the right homotopy type.  As $T$ is a
$\Sg$-inclusion prespectrum, we know that spectrification works in the
simplest possible way and that $MP$ is the homotopy colimit of the
spectra 
\[ \Sg^{-2n}\Thom(\C^n\oplus\C^n) = 
    \bigWedge_{k=-n}^n \Sg^{-2n}\Grass_{k+n}(\C^n\oplus\C^n)^\gm,
\]
where $\Grass_d(V)$ is the space of $d$-dimensional subspaces of $V$.  
It is not hard to see that the maps of the colimit system preserve
this splitting, so that $MP$ is the wedge over all $k\in\Z$ of the
spectra
\[ X_k := \hocolim_n\Sg^{-2n}\Grass_{k+n}(\C^n\oplus\C^n)^\gm. \]
This can be rewritten as
\[ X_k =
   \Sg^{2k} \hocolim_{n,m}\Sg^{-2(k+n)}\Grass_{k+n}(\C^m\oplus\C^n)^\gm.
\]
We can reindex by putting $n=i-k$ and $m=j+k$, and then pass to the
limit in $j$.  We find that
\[ X_k = \Sg^{2k}\hocolim_i\Sg^{-2i}\Grass_i(\C^\infty\oplus\C^i)^\gm. \]
It is well-known that $\Grass_i(\C^\infty\oplus\C^i)$ is a model for
$BU(i)$, and it follows that $X_k=\Sg^{2k}MU$, so
$MP=\bigWedge_k\Sg^{2k}MU$ as claimed.  We leave it to the reader to
check that the product structure is the obvious one.

All the above was done without inverting $2$.  Inverting $2$ is an
example of Bousfield localisation, and this can always be performed in
the category of strictly commutative ring spectra.

\Addresses\recd

\end{document}